\documentclass{article}

\begin{document}

\title{Lie algebra theory without algebra }
\author{S. K. Donaldson}
\maketitle
\newtheorem{thm}{Theorem}
\newtheorem{lem}{Lemma}
\newcommand{\frg}{{\bf g}}
\newcommand{\frk}{{\bf k}}
\newcommand{\frp}{{\bf p}}
\newcommand{\grad}{{\rm grad}}
\newcommand{\bR}{{\bf R}}
\newcommand{\bC}{{\bf C}}
\newcommand{\ad}{{\rm ad}}

\

\

\

\ \ \ \ \ \ \ \ \ \ \ {\it Dedicated to Professor Yu I. Manin, on his 70th. birthday.}
\

\

\section{Introduction}

  This is an entirely expository piece: the main results discussed are very well-known and the approach we take is not  really new, although
  the presentation may be somewhat different to what is in the literature.
  The author's main motivation for writing this piece comes from a feeling
  that the ideas deserve to be more widely known. 
  
  Let $\frg$ be a  Lie algebra over  $\bR$ or $\bC$. 
 A vector subspace $I\subset \frg$
  is an {\it ideal} if $[I,\frg]\subset I$. The Lie algebra  is called {\it
  simple} if it is not abelian and contains no proper ideals. A famous result
  of Cartan asserts that any simple complex Lie algebra has a compact real
  form (that is to say, the complex Lie algebra is the complexification of the Lie algebra of a compact
  group). This result underpins the theory of real Lie algebras, their maximal
  compact subgroups and the classification of symmetric spaces. 
In the standard approach, Cartan's result 
 emerges after a good deal of theory: the Theorems of Engel and Lie, Cartan's
 criterion involving the nondegeneracy of  the Killing form, root systems etc. On the other hand if one assumes this result known--by some means--then one can immediately read off much of the
standard structure theory of complex Lie groups and their representations.
Everything is reduced to the compact case (Weyl's \lq\lq unitarian trick''),
and one can proceed directly to develop the detailed theory of root systems etc. 

In \cite{kn:C}, Cartan wrote

{\it J'ai trouv\'e effectivement une telle forme pour chacun des  types de
groupes simples.
M. H. Weyl  a d\'emontr\'e ensuite l'existence de cette forme par
une raisonnement g\'en\'eral s'appliquant \`a tous les cas \`a fois. On peut se demander
si les calculs qui l'ont conduit \`a ce r\'esultat ne pourraient pas encore se
simplifier, ou plut\^ot si l'on ne pourrait pas, par une
raissonnement  a priori, d\'emontrer ce th\'eor\`eme; une telle d\'emonstration
permettrait de simplifier notablement l'exposition de la theorie des groupes simples. Je
ne suis a cet \'egard arriv\'e \`a aucun r\'esultat; j'indique simplement l'id\'ee
qui m'a guid\'e dans mes recherches infructueuses. 
}

The direct approach that Cartan outlined (in which he assumed known the nondegeneracy of the Killing form)
 was developed by Helgason (see  page 196
in \cite{kn:Hel}), and a complete proof  was accomplished by Richardson
in \cite{kn:R}. In this article we revisit these ideas and present
an almost entirely geometric proof of the result. This is essentially along the same lines
as Richardson's, so it might be asked what we can add to the story. One point
is that, guided by modern developments in Geometric Invariant Theory and
its relations with
differential geometry, we can nowadays fit this into a much more general context and hence present the proofs in a (perhaps) simpler way. Another
is that we are able to remove more of the algebraic theory; in particular,  the nondegeneracy of the Killing form. We show that the results can be deduced
 from a general principle in Riemannian geometry (Theorem 4). The arguments apply directly
 to real Lie groups, and in our exposition we will work mainly in that setting.
 In the real case the crucial concept is the following. Suppose $V$ is a
 Euclidean vector space. Then there is a transposition map $A\mapsto A^{T}$
 on the Lie algebra ${\rm End}\ V$. We say a subalgebra $\frg\subset {\rm
 End}\ V$ is {\it symmetric} with respect to the Euclidean structure if it
 is preserved by the transposition map.
 
 \begin{thm} Let $\frg$ be a simple real Lie algebra. Then there is a Euclidean
 vector space $V$, a Lie algebra embedding $\frg \subset {\rm End}(V)$, and
 a Lie group $G\subset SL(V)$ with Lie algebra $\frg$, such that
  $\frg$ is symmetric with respect to the Euclidean structure. Moreover,  any compact subgroup
 of $G$ is conjugate in $G$ to a subgroup of $G\cap SO(V)$.
 
  \end{thm}
 
 We explain in (5.1) below how to deduce the existence of the compact real form,
 in the complex case.  Theorem 1 also leads immediately to the standard
 results about real Lie algebras and symmetric spaces, as we will discuss
 further in (5.1). 
 
 The author thanks  Professors Martin Bridson, Frances Kirwan,
 Zhou Zhang and Xuhua He for comments on the earlier version of this article.

\section{More general setting}

Consider any  representation 
$$ \rho: SL(V)\rightarrow SL(W), $$
where $V,W$  are finite-dimensional real vector spaces. Let $w$ be a nonzero vector in $ W$ and
let
$G_{w}$ be the identity component of the stabiliser of $w$ in $SL(V)$.
Then we have
\begin{thm}
If $V$ is an ireducible representation of  $G_{w}$  then there is a Euclidean metric
on $V$ such that the Lie algebra of $G_{w}$ is symmetric with respect to
the Euclidean structure,  and any compact subgroup of $G_{w}$ is
conjugate in $G_{w}$ to a subgroup of $G_{w} \cap SO(V)$.
\end{thm}

Now we will show that Theorem 2 implies Theorem 1. Given a simple
real Lie algebra
$\frg$, consider the action of $SL(\frg)$ on the vector space $W$ of skew
symmetric bilinear maps from $\frg \times \frg$ to $\frg$. The Lie bracket
of $\frg$ is a point $w$ in $W$. The group $G_{w}$ is the identity
component of the group of  Lie algebra automorphisms of $\frg$, and
the Lie algebra of $G_{w}$ is the algebra ${\rm Der}(\frg)$ of derivations of $\frg$, that is, linear maps $\delta:\frg\rightarrow \frg$ with
$$  \delta [x,y]= [\delta x, y] + [x,\delta y]. $$
The adjoint action gives a Lie algebra homomorphism
$$  \ad: \frg \rightarrow {\rm Der}(\frg). $$ 
The kernel of $\ad$ is an ideal in $\frg$. This is not the whole of $\frg$
(since $\frg$ is not abelian) so it must be the zero ideal (since $\frg$ is
simple). Hence $\ad$ is injective. If $U$ is a  vector subspace of $\frg$ preserved by $G_{w}$ then any derivation $\delta$ must map $U$
to $U$. In  particular $\ad_{\xi}$ maps $U$ to $U$ for any $\xi$ in $\frg$,
so $[\frg, U] \subset U$ and $U$ is an ideal. Since $\frg$ is simple we see
that  there can
be no proper subspace preserved by $G_{w}$ and the restriction of the
representation is irreducible. By Theorem 2 there is a Euclidean metric on
$\frg$ such that ${\rm Der}(\frg)$ is preserved by transposition. Now we
want to see that in fact ${\rm Der} (\frg)=\frg$. For $\alpha\in {\rm Der}(\frg)$
and $\xi\in \frg$ we have
$$     [ ad_{\xi}, \alpha] = ad_{\alpha(\xi)}, $$
so $\frg $ is an ideal in ${\rm Der}(\frg)$. Consider the bilinear form
$$    B(\alpha_{1}, \alpha_{2}) = {\rm Tr}(\alpha_{1} \alpha_{2})$$ on
${\rm Der}(\frg)$. This is nondegenerate, since ${\rm Der}(\frg)$ is preserved
by transposition and $B(\alpha, \alpha^{T})=\vert \alpha \vert^{2}$. We have
$$   B([\alpha, \beta], \gamma)+ B(\beta, [\alpha, \gamma])=0$$
for all $\alpha,\beta, \gamma\in {\rm Der}(\frg)$. Thus the subspace
$$   \frg^{{\rm perp}} = \{ \alpha \in {\rm Der}\frg: B(\alpha, ad_{\xi})=0
\ \ {\rm for}\ {\rm all} \ \xi \in \frg\}$$
is another ideal in ${\rm Der}(\frg)$. On the other hand the map
$\alpha\mapsto - \alpha^{T}$ is an automorphism  of ${\rm Der}(\frg)$, so
$\frg^{T}$ is also an ideal in ${\rm Der}(\frg)$. Suppose that $\frg \cap
\frg^{T}\neq 0$. Then we can find a non-zero element $\alpha$ of  $\frg\cap \frg^{T}$ with $\alpha^{T} = \pm \alpha$ and then $B(\alpha, \alpha)=\pm
\vert \alpha \vert^{2}\neq 0$, so the restriction of $B$ to $\frg$ is not identically
zero. This means that $I=\frg \cap \frg^{{\rm perp}}$ is not the whole of $\frg$, but $I$ is an ideal in $\frg$ so, since $\frg$ is simple, we must
have $I=0$.  

We conclude from the above that if $\frg$ were a proper ideal in ${\rm Der}(\frg)$ there would be another proper ideal $J$ in ${\rm Der}(\frg)$
such that $J\cap \frg=0$. (We take $J$ to be either $\frg^{T}$ or $\frg^{{\rm
perp}}$.) But then for $\alpha \in J$ we have $[\alpha, \frg]=0$, but this
means that $\alpha$ acts trivially on $\frg$, which gives a contradiction.

Finally, the statement about compact subgroups in Theorem 1 follows immediately
from that in Theorem 2.

\

(The argument corresponding to the above in the complex case (see (5.1)) is more transparent. )

 \section{Lengths of vectors}
 We will now begin the proof of Theorem 2. The idea is to find a metric
 by minimising the associated norm of the vector $w$. In the Lie algebra
 situation, which we are primarily concerned with here, this is in essence
 the approach suggested by Cartan and carried through by Richardson. In the
 general situation considered in Theorem 2 the ideas have been studied and
 applied extensively over the last quarter century or so, following the work
 of Kempf-Ness \cite{kn:KN}, Ness \cite{kn:N} and Kirwan \cite{kn:K}. Most of the literature
 is cast in the setting of complex representations. The real case has
 been studied by Richardson and Slodowy \cite{kn:RS} and Marian  \cite{kn:M} and works in just
 the same way.

 Recall that we have a representation $\rho$ of $SL(V)$ in $SL(W)$, where
 $V$ and $W$ are real vector spaces, a fixed vector $w\in W$ and we define
 $G_{w}$ to the stabiliser of $w$ in $SL(V)$. Suppose we also have some compact
 subgroup (which could be trivial) $K_{0}\subset G_{w}$. We fix any Euclidean metric $\vert \
 \vert_{1}$ on $V$ which is preserved by $K_{0}$. Now it is standard that
 we can choose a Euclidean metric $\vert \ \vert_{W}$ on $W$ which is invariant under the restriction
 of $\rho$ to $SO(V)$. We want to choose this metric $\vert \ \vert_{W}$
 with the further property that the derivative $d\rho$ intertwines transposition
 in ${\rm End} V$ (defined by $\vert \ \vert_{1}$) and transposition in ${\rm
 End} W$ (defined by $\vert \ \vert_{W}$); that is to say
 $$  d\rho(\xi^{T}) = \left( d\rho(\xi) \right)^{T}. $$
 To see that this is possible we can argue as follows. We  complexify
 the representation to get $\rho_{\bC}: SL(V\otimes \bC)\rightarrow SL(W\otimes
 C)$. Then the compact group generated by the action of $SU(V\otimes \bC)$
 and complex conjugation acts on $W\otimes \bC$ and we can choose a Hermitian
 metric on $W\otimes \bC$ whose norm function is invariant under this group.
 Invariance under complex conjugation means that this Hermitian metric is
 induced from a Euclidean metric on $W$. Then the fact that $\rho_{\bC}$
 maps $SU(V\otimes \bC)$ to $SU(W\otimes \bC)$ implies that 
 $d\rho$ has the property desired. (The author is grateful to Professors
 He and Zhang for pointing out the need for this argument. In our main application,
 to Theorem 1, 
 the standard metric on $W$ already has the desired property.)

Now define a function $\tilde{F}$ on $SL(V)$ by
 $$  \tilde{F}(g)= \vert g(w)\vert_{W}^{2}.$$
 For $u\in SO(V)$ and $\gamma \in G_{w}$ we have $$ \tilde{F}(ug \gamma)= \vert u g \gamma (w)\vert_{W}^{2} = \vert
 u g(w)\vert_{W}^{2} = \tilde{F}(g) $$
 So $\tilde{F}$ induces a function $F$ on the quotient space ${\cal H}=SL(V)/SO(V)$,
 invariant under the natural action of $G_{w} \subset SL(V)$. We can think
 about this in another, equivalent, way. We identify ${\cal H}$ with the Euclidean
 metrics on $V$ of a fixed determinant. Since $\rho:SL(V)\rightarrow SL(W)$
 maps $SO(V)$ to $SO(W)$ it induces a map from $SL(V)/SO(V)$ to $SL(W)/SO(W)$
 and so a metric on $V$ with the same determinant as $\vert \ \vert_{1}$ induces a metric
 on $W$. Then  $F$ is given by the square of the induced norm of the fixed vector
 $w$. Explicitly, the identification of $SL(V)/SO(V)$ with  metrics is given
 by $[g]\mapsto \vert \ \vert_{g}$ where
 $$   \vert v \vert_{g}^{2} = \vert gv \vert_{1}^{2} = \langle v, g^{T} g
 v \rangle_{1}. $$
 
 This function $F$  has two crucial, and well-known, properties, which we  state in the following Lemmas
 
  \begin{lem}
  Suppose $F$ has a critical point at $H\in {\cal H}$. Then the Lie algebra
  of the stabiliser $G_{w}$ is symmetric with respect to the Euclidean structure
  $H$ on $V$.
  \end{lem}
  To prove this, there is no loss in supposing that $H$ is the original metric
  $\vert\ \vert_{1}$. (For we can replace $w$ by $g w$ for any $g\in SL(V)$.)
  The fact that $\rho$ maps $SO(V)$ to $SO(W)$ implies that its derivative
  takes transposition in ${\rm End} V$ (defined by $\vert \ \vert_{1}$) to
  transposition in ${\rm End} W$ defined by $\vert \ \vert_{W}$. The
  condition for $\tilde{V}$ to be stationary is that
  $$  \langle d\rho(\xi) w, w\rangle_{W} =0 $$
  for all $\xi$ in the Lie algebra of $SL(V)$. In particular consider elements
  of the form $\xi=[\eta, \eta^{T}]$ and write $A=d\rho(\eta)$. Then we have
  $$ 0=\langle d\rho [\eta, \eta^{T} ] w, w \rangle_{W}= \langle [A,A^{T}] w,w\rangle_{W}=
  \vert A^{T} w \vert_{W}^{2} - \vert A w \vert_{W}^{2}. $$
  By definition $\eta$ lies in the Lie algebra of $G_{w}$ if and only if
  $A w=0$. By the identity above, this occurs if and only if $A^{T} w=0$,
  which is just when $\eta^{T}$ lies in the Lie algebra of $G_{w}$.  
  
  \
  
  For the second property of the function we need to recall the standard
  notion of geodesics in ${\cal H}$.  We can identify ${\cal H}$ with the positive
  definite  symmetric elements of $SL(V)$, with the quotient map
  $SL(V)\rightarrow {\cal H}$ given by $g\mapsto g^{T} g$. Then the geodesics
  in ${\cal H}$ are paths of the form   
 
 \begin{equation} \gamma(t)= g^{T} \exp(St) g, \label{eq:geodesic} \end{equation}
 where $g$ and $S$ are fixed, with $g\in SL(V)$ and $S$ a trace-free
 endomorphism which is symmetric with respect to  $\vert \ \vert_{1}$. Another
 way of expressing this is that a geodesics through any point $H\in {\cal
 H}$ is the orbit of $H$ under a $1$-parameter subgroup $e(t)$ in $ SL(V)$ where
 $ e(t) = \exp(\sigma t)$ with $\sigma$ a symmetric endomorphism with
 respect to the metric $H$. 
 \begin{lem}
 \begin{enumerate}
 \item For any geodesic $\gamma$ the function $F\circ\gamma$ is convex i.e.
 $$ \frac{d^{2}}{dt^{2}} F(\gamma(t)) \geq 0. $$
 \item If $F$ achieves its minimum in ${\cal H}$ then $G_{w}$ acts transitively on the set of minima.
 \end{enumerate}
 \end{lem} 

 To prove the first part, note that, replacing $w$ by $gw$, we can reduce to considering a geodesic through the base point $[1]\in {\cal H}$, so of
the form $\exp(St)$ where $S$ is symmetric with respect to $\vert \ \vert_{1}$.
   Now the derivative $d\rho$ maps the symmetric endomorphism $S$ to a symmetric
 endomorphism $A\in {\rm End}(W)$. We can choose an orthonormal basis in $W$
 so that $A$ is diagonal, with eigenvalues $\lambda_{i}$ say. Then
 if $w$ has coordinates $w_{i}$ in this basis we have
 $$  F(\exp(St))= \tilde{F}(\exp(St/2))= \sum \vert \exp(\lambda_{i} t/2) w_{i}\vert_{W}^{2}= \sum
 \vert w_{i} \vert^{2} \exp( \lambda_{i} t),$$
 and this is obviously a convex function of $t$.
 
 To prove the second part note that, in the above, the function 
 $F(\exp(St)$ is either strictly convex or constant, and the latter only
 occurs when  $\lambda_{i}=0$ for each index $i$ such that $w_{i}\neq 0$,
 which is the same as saying that $\exp(St) w=w$ for all $t$, or that
 the $1$-parameter subgroup $\exp(St)$ lies in $G_{w}$. More generally if
 we write a geodesic through a point $H$ as the orbit of $H$ under a $1$-parameter
 subgroup $e(t)$ in $SL(V)$ then the function is constant if and only if
 the $1$-parameter subgroup lies in $G_{w}$. Suppose that
 $H_{1}, H_{2}$ are two points in ${\cal H}$ where $F$ is minimal. Then
 $F$ must be constant on the geodesic between $H_{1}, H_{2}$. 
 Thus $H_{2}$ lies in the orbit of $H_{1}$ under a $1$-parameter subgroup
 in $G_{w}$. So $G_{w}$ acts transitively on the set of minima.
 
 \

 We now turn back to the proof of Theorem 2. Suppose that the convex function
 $F$ on ${\cal H}$ achieves a minimum at $H_{1}\in {\cal H}$. Then by Lemma 1 the Lie algebra
 of $G_{w}$ is symmetric with respect to the Euclidean structure $H_{1}$ on $V$. It only remains
 to see that the compact subgroup $K_{0}$  of $G_{w}$ is conjugate to a subgroup of the orthogonal
 group for this Euclidean structure.
 For each $H\in {\cal H}$ we have a corresponding special orthogonal group $SO(H,V)\subset SL(V)$. For $g\in SL(V)$ the groups $SO(H,V), SO(g(H),V)$
are conjugate by $g$ in $SL(V)$.
 Recall that we chose the metric $\vert
 \ \vert_{1}$ to be $K_{0}$ invariant. This means that $K_{0}$ fixes the
 base point $[1]$ in ${\cal H}$. Suppose we can find a point $H_{0}$ in ${\cal
 H}$ which minimises $F$ and which is also $K_{0}$-invariant. Then 
 $K_{0}$ is contained in $SO(H_{0},V)$. But by the second part of Lemma
 2 there is a $\gamma\in G_{w}$ such that $\gamma(H_{0})=H_{1}$. Thus conjugation
 by $\gamma$ takes $SO(H_{0}, V)$ to $SO(H_{1}, V)$ and takes $K_{0}$ to
 a subgroup on $SO(H_{1}, V)$, as required.
 
 To sum up, Theorem 2 will be proved if we can establish the following result.
   
 \begin{thm}
 Let $F$ be a convex function on ${\cal H}$, invariant under a group
 $G_{w} \subset SL(V)$. Let $K_{0}$ be a compact subgroup of $G_{w}$ and let $[1]\in {\cal H}$ be fixed by $K_{0}$. Then if $V$ is an ireducible representation of
 $G_{w}$ there is a point $H_{0}\in {\cal H}$ where $F$ achieves its minimum
 and which is fixed by $K_{0}$.
 \end{thm}
 
  (Notice that the hypothesis here that there is a point $ [1]\in {\cal H}$ fixed
  by $K_{0}$ is actually redundant, since any compact subgroup of $SL(V)$
  fixes some metric.)

\section{Riemannian geometry argument}

In this section we will see that Theorem 3 is a particular case of a more
general result in Riemannian geometry. Let $M$ be a complete Riemannian manifold,
so for each point $p\in M$ we have a surjective exponential map
$$  \exp_{p}: TM_{p} \rightarrow M. $$
We suppose $M$ has the following property

\

{\bf Property (*)}

{\it For each point $p$ in $M$ the exponential map $\exp_{p}$ is distance-increasing}

\

Readers with some background in Riemannian geometry will know that it is
equivalent to say that $M$ is simply connected with nonpositive sectional
curvature, but we do not need to assume knowledge of these matters. The crucial
background we need to know is

\

{\bf Fact}

{\it There is a metric on ${\cal H}= SL(V)/SO(V)$ for which the action
of $SL(V)$ is isometric, with the geodesics described in (1) above
and having Property (*).}

 \
 
 This Riemannian metric
on ${\cal H}$ can be given by the formula
$$   \Vert \delta H \Vert^{2}_{H} = {\rm Tr} \left(\delta H H^{-1}\right)^{2}.$$
  The distance-increasing property can be deduced from the fact
that ${\cal H}$ has non-positive curvature and standard comparison results
for Jacobi fields. For completeness, we give a self-contained proof of the Fact in
the Appendix. 

\

The piece of theory we need to recall in order to state our Theorem is the notion of the \lq\lq sphere
at infinity'' associated to a manifold $M$ with Property (*).
This will be familiar in the prototype cases of Euclidean space and hyperbolic
space. In general, for $x\in M$ write $S_{x}$ for the unit sphere in the
tangent space $TM_{x}$ and define 
 $$ \Theta_{x}: M\setminus \{x\} \rightarrow S_{x}$$
   by $$\Theta_{x}(z)= \frac{\exp_{x}^{-1}(z)}{\vert \exp_{x}^{-1}(z)\vert}.
   $$If $y$ is another point in $M$ and $R$ is greater than the distance
   $d= d(x,y)$ we define
   $$  F_{R,x,y}: S_{y}\rightarrow S_{x}$$
   by $$ F_{R,x,y}(\nu)= \Theta_{x} \exp_{y}(R \nu). $$
   \begin{lem}
   For fixed $x,y,\nu$ the norm of the derivative of $F_{R,x,y,\nu}$ with respect
   to $R$ is bounded by
   $$   \vert \frac{\partial}{\partial R}F_{R,x,y}(\nu) \vert \leq \frac{d}{R
   (R-d)}.$$
   
   \end{lem}
   
   Let $\gamma$ be the geodesic $\gamma(t)= \exp_{y}(t \nu)$, let $w$ be the
   point $\gamma(R)$ and let $\sigma$ be the geodesic from $x$ to $w$. The distance-increasing property of $\exp_{x}$ implies that the norm of the derivative
appearing in the statement is bounded by $d(x,w)^{-1}$ times the component
of $\gamma'(R)$ orthogonal to the tangent vector of $\sigma$ at $w$. Thus
$$ \vert \frac{\partial }{\partial R}{F_{R,x,y}(\nu)}\vert \leq \frac{\sin \phi}{
d(x,w)}, $$
where $\phi$ is the angle between the geodesics $\gamma, \sigma$ at $w$.
By the triangle inequality $d(x,w)\geq R-d$. In a Euclidean triangle with
side lengths $d,R$ the angle opposite to the side of length $d$ is at most $\sin^{-1}(d/R)$.
It follows from the distance-increasing property of $\exp_{z}$ that $\sin
\phi\leq d/R$. Thus
$$ \frac{\sin \phi}{d(x,w)} \leq \frac{d}{R(R-d)}, $$
as required.

   Since the integral 
   of the function $1/R(R-d)$, with respect to $R$, from $R=2d$ (say) to $ R=\infty$, is finite,
   it follows
  from the Lemma that $F_{R,x,y}$ converges uniformly as $R\rightarrow
 \infty$ to a continuous map $F_{x,y}:S_{y}\rightarrow S_{x}$, and obviously
 $F_{x,x}$ is the identity. Let $z$ be
 another point in $M$ and $\nu$ be a unit tangent vector at $z$. Then we
 have an identity, which follows immediately from the definitions,
 $$   F_{R, x,z}(\nu)= F_{R',x, y} \circ F_{R,y,z} (\nu), $$
 where $R'= d(y, \exp_{z}(R \nu))$. Since, by the triangle inequality again,
 $$    R'\geq R- d(y,z), $$
 we can take the limit as $R\rightarrow \infty$ to obtain
 $$  F_{x,z} = F_{x,y} \circ F_{y,z} : S_{z} \rightarrow S_{x}.$$
 In particular, $F_{y,x}$ is inverse to $F_{x,y}$ so the maps $F_{x,y}$ give
 a compatible family of homeomorphisms between spheres in the tangent spaces.
 We define the {\it sphere at infinity} $ S_{\infty}(M)$ to be the quotient
 of the unit sphere bundle of $M$ by these homeomorphisms, with the topology
 induced by the identification with $S_{x_{0}}$ for any fixed base point $x_{0}$.

 Now suppose that a topological group $\Gamma$ acts by isometries on
   $M$. Then $\Gamma $ acts on $S_{\infty}(M)$, as a set.       
 Explicitly, if we fix a base point $x_{0}$ and identify the sphere at infinity
 with $S_{x_{0}}$,, the action  of a group element $g\in \Gamma$ is given by
 $$  g(\nu)= \lim_{R\rightarrow \infty} \Theta_{x_{0}} g (\exp_{x_{0}} R
 \nu). $$
 Write the action as $$ A: \Gamma\times S_{x_{0}} \rightarrow S_{x_{0}}.
 $$
 Given a compact set $P\subset \Gamma$ we can define
 $$  A_{R}: P \times S_{x_{0}} \rightarrow S_{x_{0}}, $$
 for sufficiently large $R$, by
 $$  A_{R}(g,\nu)= \Theta_{x_{0}} g(\exp_{x_{0}} R \nu). $$
 Since $g(\exp_{x_{0}} R \nu)= \exp_{g(x_{0})}(R g_{*} \nu)$ the maps $A_{R}$ converge
 uniformly as $R\rightarrow \infty$ to the restriction of $A$ to $P\times
 S_{x_{0}}$. It follows that the action $A$ is continuous.
 With these preliminaries in place we can state our main technical result.

 \begin{thm} Suppose that the Riemannian manifold $M$ has Property (*). Suppose
 that $\Gamma$ acts by isometries on $M$ and $F$ is a convex $\Gamma$-invariant
 function on $M$. Then either there is a fixed point for the action of $\Gamma$
 on $S_{\infty}(M)$ or the function $F$ attains its minimum in $M$. Moreover,
 in the second case, if $K_{0}$ is a subgroup of $\Gamma$ which fixes a point
 $x\in M$, then there is a point  $x'\in M$ where $F$ attains its minimum in
 $M$ and with $x'$ fixed by $K_{0}$.
 \end{thm}

 Return now to our example ${\cal H}$. The tangent space at the identity
 matrix $[1]$ is the set of trace-free symmetric matrices.
  We define a {\it weighted flag}
 $({\cal F}, \underline{\mu})$ to be a strictly increasing sequence of vector
 subspaces
 $$  0=F_{0}\subset  F_{1} \subset F_{2} \dots \subset F_{r} = V$$
 with associated weights $\mu_{1}>\mu_{2} \dots >\mu_{r}$, subject to the
 conditions
 $$ \sum n_{i} \mu_{i} = 0, \sum n_{i} \mu_{i}^{2} =1, $$
 where $n_{i}= {\rm dim} F_{i}/F_{i-1}$. If $S$ is a trace-free symmetric
 endomorphism
 with ${\rm Tr}\ S^{2}=1$ then we associate a weighted flag to $S$
 as follows.
 We take $\mu_{i}$ to  be the eigenvalues of $S$, with eigenspaces $E_{i}$,
 and
 form a flag with
 $$ F_{1}= E_{1}\ ,\ F_{2} = E_{1} \oplus E_{2} , \dots. $$
 It is clear then that the unit sphere $S_{[1]}$ in the tangent space of ${\cal
 H}$ at $[1]$ can be identified with the set of all weighted flags. Now
 there is an obvious action of $SL(V)$ on the set of weighted flags and
 we have:
 \begin{lem}
 The action of $SL(V)$ on the sphere at infinity in ${\cal H}$
 coincides with the obvious action under the identifications above.
 \end{lem}
 This is clearly true for the subgroup $SO(V)$. We use the fact that given
 any weighted flag $({\cal F}, \underline{\mu})$ and $g\in SL(V)$ we can write $g=u h$ where
 $u\in SO(V)$ and $h$ preserves $\cal F$. (This is a consequence of the obvious
 fact that $SO(V)$ acts transitively on the set of flags of a given type.)
 Thus it suffices to show that
  such $h$ fix the point $S$ in the unit sphere corresponding to $({\cal F},\underline{\mu})$
 in the differential-geometric action. By the $SO(V)$ invariance of the set-up
 we can choose a basis so that  ${\cal F}$ is the standard flag
 $$  0 \subset \bR^{n_{1}} \subset \bR^{n_{1}}\oplus \bR^{n_{2}} \dots \subset
 \bR^{n}. $$
 Then $S$ is the diagonal matrix with diagonal entries $\mu_{1},\dots, \mu_{r}$, repeated according to the multiplicities $n_{1}, \dots, n_{r}$. The matrix $h$ is upper triangular
 in blocks with respect to the flag. Now consider, for a large real parameter $R$
 the matrix
 $$  M_{R} = \exp(-\frac{R S}{2}) h \exp (\frac{RS}{2}).  $$
Consider a block $h_{ij}$ of $h$. The corresponding block of $M_{R}$ is
$$  \left(M_{R}\right)_{ij}= e^{R(\mu_{i}-\mu_{j})/2} h_{ij}. $$
Since  $h$ is upper-triangular in blocks and  the $\mu_{i}$ are increasing,
we see that $M_{R}$ has a limit as $R$ tends to infinity, given by the diagonal
blocks in $h$. Since these diagonal blocks are invertible the limit of $M(R)$
is invertible,  hence
$$\delta_{R}= {\rm Tr}\left( \log ( M_{R} M_{R}^{*})\right)^{2}$$ is a bounded
function of $R$.  But $\delta_{R}^{1/2}$ is the distance in ${\cal H}_{n}$ between $\exp(RS)$
and $h \exp(RS) h^{T}$. It follows from the comparison argument, as before,
that the angle between $\Theta_{[1]} (h \exp(RS) h^{T})$ and $S$ tends to zero
as $R\rightarrow \infty$, hence $h$ fixes $S$ in the differential geometric
action.

\

 Now Theorem 3 is an immediate consequence of Theorem 4 and Lemma 4, since
 if $G_{w}$ fixes a point on the sphere at infinity in ${\cal H}$ it fixes
 a flag, hence some non-trivial subspace of $V$, and $V$ is reducible as a
 representation of $G_{w}$. 
 
 \pagebreak
 
 {\bf Remarks}
 \begin{itemize}
 \item
 The advantage of this approach is that Theorem 4 seems quite accessible
 to geometric intuition. For example it is obviously true in the case when
 $M$ is hyperbolic space, taking the ball model, and we suppose that $F$ extends
 continuously to the boundary of the ball. For then $F$ attains its minimum
 on the closed ball and if there are no minimising points in the interior
 the minimiser on the boundary must be unique (since there is a geodesic
 asymptotic to any two given points in the boundary). 
 
 \item The author has not found Theorem 4  in the
 literature, but it does not seem likely that it is new. There are very similar
 results in \cite{kn:BOn} for example. The author has been told by
 Martin Bridson that a more general result of this nature holds, in the context
 of proper CAT(0) spaces. The proof of this more general result follows in
 an obvious way from Lemma 8.26 of \cite{kn:BH} (see also Corollary 8.20
 in that reference).
 
 \item The hypothesis on the existence of a fixed point $x$ for $K_{0}$ in
 the statement of Theorem 4 is redundant, since any compact group acting
 on a manifold with Property (*) has a fixed point, by a theorem of Cartan
 (see the remarks at the end of Section 3 above, and at the end of (5.2) below).
 However we do not need to use this.  
 \end{itemize}
 
 \
 
 We now prove Theorem 4. We begin by disposing of the statement involving
 the compact group $K_{0}$. Suppose that $F$ attains its minimum somewhere
 in $M$. Then, by convexity, the minimum set is a totally geodesic submanifold
 $\Sigma \subset M$. The action of $K_{0}$ preserves $\Sigma$, since $F$ is
 $\Gamma$-invariant and $K_{0}$ is contained in $\Gamma$. Let
 $x'$ be a point in $\Sigma$ which minimises the distance to the $K_{0}$
 fixed-point $x$. Then if
 $x''$ is any other point in $\Sigma$ the geodesic segment from $x'$ to $x''$
 lies in $\Sigma$ and is orthogonal to the geodesic from $x$ to $x'$ at $x'$.
 By the distance-increasing property of the exponential map at $x'$ it follows
 that the distance from $x$ to $x''$ is strictly greater than the distance
 from $x$ to $x'$. Thus the distance-minimising point $x'$ is unique, hence fixed
 by $K_{0}$.
  
  \
  
 To prove the main statement in  Theorem 4  we use the following Lemma.
 \begin{lem}
 Suppose that $M$ has Property (*) and $N$ is any set  of isometries of $M$.
If  there is a sequence $x_{i}$ in $M$ with $d(x_{0}, x_{i})\rightarrow \infty$
 and for each $g\in N$ there is a $C_{g}$ with  $d(x_{i}, g x_{i})\leq C_{g}$
 for all $i$, then there is a point in $S_{\infty}(M)$
 fixed by $N$.
 
 \end{lem}

 Set $R_{i}= d(x_{0}, x_{i})$ and $\nu_{i}= \Theta_{x_{0}}(x_{i}) \in S_{x_{0}}$.
 By the compactness of this sphere we may suppose, after perhaps taking a
 subsequence, that the $\nu_{i}$ converge as $i$ tends to infinity to some
 $\nu\in S_{x_{0}}$.
 Then for each $g\in N$ we have, from the definitions,
 $$  A_{R_{i}}(g, \nu_{i})= \Theta_{x_{0}} ( g x_{i}). $$
 Fixing $g$, let $\phi_{i}$ be the angle between the unit tangent vectors
 $\nu_{i}=\Theta_{x_{0}}(x_{i}) $ and $\Theta_{x_{0}}( g x_{i})$. The distance
 increasing property implies, as in Lemma 3, that
 $$  \sin \phi_{i} \leq C_{g}/ R_{i}. $$
 The angle $\phi_{i}$ can be regarded as the distance ${\rm dist}(\ ,\ )$ between the points
 $\nu_{i}$ and $A_{R_{i}}(g, \nu_{i})$ in the sphere $S_{x_{0}}$. In other
 words we have
 $$ {\rm dist}\ ((\nu_{i}, A_{R_{i}}(g, \nu_{i})) \leq \sin^{-1}( \frac{C_{g}}{R_{i}})
 .$$
  Now take the limit as $i\rightarrow \infty$: we see that
  $ {\rm dist}\ (\nu, A(g, \nu))=0$, which is to say that $\nu$ is fixed by $g$.

 \
 
 To prove Theorem 4, consider the gradient vector field $\grad
 \ F$ of the function $F$, and the associated flow
 $$   \frac{dx}{dt}  = -\grad\ F_{x}$$
 on $M$. By the standard theory, given any initial point there is a solution
 $x(t)$ defined for some time interval $(-T,T)$.
 \begin{lem}
 If $x(t)$ and $y(t)$ are two solutions of the gradient flow equation,
 for $t\in (-T,T)$, then
  $d(x(t), y(t))$ is a non-increasing function of $t$.
 \end{lem}
 If $x(t)$ and $y(t)$ coincide for some $t$ then they must do so for all
 $t$, by uniqueness to the solution of the flow equation, and in that case
 the result is certainly true. If $x(t)$ and $y(t)$ are always different then
 the function $D(t)= d(x(t), y(t))$ is smooth: we compute the derivative at some
 fixed $t_{0}$. Let $\gamma(s)$ be the geodesic from $x(t_{0})= \gamma(0)$ to $y(t_{0})= \gamma(D)$. Clearly
$$ D'(t_{0})= \langle \grad F_{x(t)}, \gamma'(0)\rangle - \langle \grad F_{y(t)},
\gamma'(D)\rangle. $$
But $\langle \grad F_{\gamma(s)}, \gamma'(s)\rangle$ is the derivative of the function
$F \circ\gamma(s)$, which is nondecreasing in $s$ by the convexity hypothesis, so
$D'(t_{0}) \geq 0$, as required.
  
A first consequence of this Lemma---applied to $y(t)=x(t+\delta)$ and taking
the limit as $\delta \rightarrow 0$--- is that
the velocity $\vert \frac{dx}{dt}\vert $of a gradient path is decreasing.
Thus for finite positive time $x(t)$ stays in an {\it a priori} determined compact
subset of $M$ (since this manifold is complete). It follows that the flow
is actually defined for all positive time, for any initial condition.  Consider
an arbitrary initial point $x_{0}$ and let $x(t)$ be this gradient path,
for $t\geq 0$. If there is a sequence $t_{i}\rightarrow \infty$ such that
$x(t_{i})$ is bounded then, taking a subsequence, we can suppose that
$x(t_{i})$ converges and it follows in a standard way that the limit is a
minimum of $F$. If there is no such minimum then we can take a sequence such that  $x_{i}= x(t_{i})$ tends to infinity. Suppose that $g$ is in $ \Gamma$, so the action of $g$ on $M$ preserves
$F$ and the metric. Then $y(t)= g(x(t))$ is a gradient path with initial value $g(x_{0})$
and  $d(x_{i}, g x_{i})\leq C_{g} = d(x_{0}, g x_{0})$. Then, by Lemma 5, there is a fixed point for the action of $\Gamma$ on $S_{\infty}(M)$.

\

There is an alternative argument which is perhaps more elementary, although
takes more space to write down in detail. With a fixed base point $x_{0}$ choose $c$
with $\inf_{M} F < c<
F(x_{0})$ and let $\Sigma_{c}$ be the  hypersurface $F^{-1}(c)$. Let
$z_{c} \in \Sigma_{c}$ be a point which minimises the distance to $x_{0}$, so
that $x_{0}$ lies on a geodesic $\gamma$ from $z_{c}$ normal to $\Sigma_{c}$. The convexity
of $F$ implies that the second fundamental form of $\Sigma_{c}$ at $z_{c}$ is
positive with respect to the normal given by the geodesic from $z_{c}$ to
$x_{0}$. A standard comparison argument in \lq\lq Fermi coordinates'' shows
that the exponential map on the normal bundle of $\Sigma_{c}$ is distance increasing
on the side towards $x_{0}$. In particular, let $w$ be another point in $\Sigma_{c}$
and $y = \exp (R \xi)$ where $\xi$ is the unit normal to $\Sigma_{c}$ at $w$  pointing in the direction of increasing
$F$ and $R=d(x_{0}, z_{c})$. Then we have $d(z_{c}, w) \leq d(x_{0}, y)$. Now suppose $g$ is in $\Gamma$.
Then $g$ preserves $\Sigma_{c}$ and if we take $w=g(z_{c})$ above we have $y=g(x_{0})$.
So we conclude from this comparison argument that $d(z_{c}, g(z_{c})) \leq
d(x_{0}, g x_{0})$. Now take a sequence $c_{i}$ decreasing to $\inf F$ (which
could be finite or infinite). We get a sequence $x_{i} = z_{c_{i}}$ of points
in $M$. If $(x_{i})$ contains  a bounded subsequence then we readily deduce that there
is a minimum of $F$. If   $x_{i}$ tends to infinity we get a sequence to
which we can apply Lemma 5, since $d(x_{i}, g x_{i}) \leq C_{g} = d(x_{0},
g x_{0})$.      

\section{Discussion}
 
 \subsection{Consequences of Theorem 1}
 \begin{itemize}
 \item We start with a simple Lie algebra $\frg$ and use Theorem 1 to obtain an
 embedding $\frg\subset {\rm End}(V)$, for a Euclidean space $V$, with $\frg$
 preserved by the transposition map. We also have a corresponding Lie group $G \subset SL(V)$.
 We write $K$ for the identity component of $G\cap SO(V)$. It follows immediately
  from
 Theorem 1 that $K$ is a maximal compact connected subgroup of $G$, and any
 maximal compact connected subgroup is conjugate to $K$.
 
 \item The involution $\alpha \mapsto - \alpha^{T}$ on ${\rm End}(V)$ induces
 a Cartan involution of $\frg$ so we have an eigenspace decomposition
 $$   \frg = \frk \oplus \frp $$
 with $\frk ={\rm Lie}(K)$ and
 \begin{equation}  [\frk,\frk]\subset \frk\  ,\  [\frk,\frp]\subset \frp\
 ,\  [\frp, \frp] \subset
 \frk. \label{eq:inclusions} \end{equation}
 Notice that $\frk$ is non-trivial, for otherwise $\frg$ would be abelian.
 \item Consider the bilinear form $B(\alpha, \beta)= {\rm Tr}(\alpha \beta)$ on
 $\frg$. Clearly this is positive definite on $\frp$, negative definite on $\frk$
 and the two spaces are $B$-orthogonal. Thus $B$ is nondegenerate. The Killing
 form $\hat{B}$ of $\frg$ is negative-definite on $\frk$ (since the restriction of
 the adjoint action to $K$ preserves some metric and $\frk$ is not an ideal).
 So the Killing form is not identically zero and must be a positive multiple
 of $B$ (otherwise the relative eigenspaces would be proper ideals). In fact
 we do not really need this step, since in our proof of Theorem 1 the vector
 space $W$ is $\frg$ itself, and $B$ is trivially equal to the Killing form. 
  \item Either $\frp$ is trivial, in which case $G$ is itself compact, or
  there is a nontrivial Riemannian symmetric space of negative type $M_{\frg}^{-}=G/K$
  associated to $\frg$. This can be described rather explicitly. Let us now
  fix on  the specific representation $\frg \subset {\rm End}(\frg)$ used in
  the proof of Theorem 1. Say a Euclidean metric on $\frg$ is \lq\lq optimal''
  if the adjoint embedding is symmetric with respect to the metric, as in Theorem
  1. (It is easy to see that the optimal metrics are exactly those which minimise the
  norm of the bracket, among all metrics of a given determinant.) Then
  $M^{-}_{\frg}$ can be identified with the set of optimal metrics, a totally
  geodesic submanifold of ${\cal H}= SL(\frg)/SO(\frg)$.
  
  \item So far we have worked exclusively in the real setting. We will now
  see how to derive the existence of compact real forms of a simple complex Lie
  algebra.    
 
 \begin{lem} If $\frg$ is a simple complex Lie algebra then it is also simple
 when regarded as a real Lie algebra. \end{lem}
 
 To see this, suppose that $A\subset \frg$ is a proper real ideal: a real vector
 subspace with $[A,\frg]\subset A$. By complex linearity of the bracket
 $A\cap iA$ is a complex ideal, so
 we must have $A\cap iA=0$. But then since $i\frg  = \frg $ we have
 $[A,\frg]= [A, i\frg]=i [A,\frg]\subset iA$, so $[A,\frg]=0$. But  $A+ i A$ is another
 complex ideal, so we must have $\frg= i A \oplus A$ and $\frg$ is Abelian.
 
 Next we have
 \begin{lem}
 Let $\frg$ be a simple complex Lie algebra and let  $\frg= {\rm Lie}(G)\subset
 {\rm End} V $ be an embedding provided by Theorem 1, regarding $\frg$ as a
 real Lie algebra. Then $\frg$ is the complexification of the Lie algebra
 of the compact group  $K=G\cap SO(V)$.
 \end{lem}
The inclusions (2) imply that

 $$  I= (\frp \cap i \frk)+ (\frk\cap i \frp), $$
 is a complex ideal in $\frg$, so either $I=\frg$ or $I=0$. In the first
 case we have $i\frk=\frp$ and $\frg$ is the complexification of $\frk$,
 as required. So we have to rule out the second case. If this were to hold
 we have $\frk \cap i \frp=0$ so $\frg= \frk \oplus (i \frp)$. Then
 $[i\frp, i\frp] \subset \frk $ so the map $\sigma$ on $\frg$ given by multiplication
 by $1$ on $\frk$ and by $-1$ on $i\frp$ is another involution of $\frg$,
 regarded as a real Lie algebra. Now let  $\hat{B}_{\bC}$ be
 the Killing form regarded as a complex Lie algebra. So $\hat{B} = 2 {\rm
 Re}\hat{B}_{\bC}$. The fact that $\sigma$ is an involution of $\frg$ means that
 $\hat{B}(\frk, i\frp)=0$. But we know that $\frp$ is the orthogonal complement
 of $\frk$ with respect to $B$ and $\hat{B}$, so we must have $i\frp = \frp$. But
 $\hat{B}$ is positive definite on $\frp$ while $\hat{B}(i\alpha, i\alpha) =
 2{\rm Re} \hat{B}_{\bC}(i\alpha, i\alpha) = - \hat{B}(\alpha,\alpha)$ so
 $\frp\cap i\frp=0$. This means that $\frp=0$ and $\frg=\frk$ which is clearly
 impossible (by the same argument with the Killing form).
 
 \item The argument above probably obscures the picture. If one is interested
 in  the
 complex situation it is much clearer to redo the whole proof in this setting,
 working with Hermitian metrics on complex representation spaces. The proof
 goes through essentially word-for-word, using the fact that the standard
 metric on $SL(n,\bC)/SU(n)$
 has Property (*).Then one can deduce the real case from the complex case
 rather than the other way around, as we have done above.
 
 \item Returning to the case of a simple real Lie algebra $\frg$, which is
 not the Lie algebra of a compact group,  we can
 also give an explicit description of the symmetric space $M^{+}_{\frg}$
 of positive type dual to $M^{-}_{\frg}$. Fix an optimal metric on $\frg$
 and extend it to a Hermitian metric $H$ on $\frg \otimes \bC$. Then
 $M^{+}_{\frg}$ is the set of real forms $\frg'\subset \frg\otimes \bC$ which
 are conjugate by $G^{c}$ to $\frg$ and such that the restriction of ${\rm
 Re}\ H$
 to $\frg'$ is an optimal metric on $\frg'$. This is a totally geodesic submanifold
 of $SU(\frg\otimes \bC)/SO(\frg)$. 
 
 \end{itemize}
              
 \subsection{Comparison with other approaches}
 
 The approach we have used, minimising the norm of the Lie bracket, is essentially
 the same as that suggested by Cartan, and carried through by Richardson,
 with the difference that we do not assume known that the Killing form is
 nondegenerate so we operate with a special linear group rather than an orthogonal
 group. The crucial problem is to show that the minimum is attained when the
 Lie algebra is simple.  This can be attacked by considering points in the
 closure of the relevant orbit for the action on the projectivized space.
 Richardson gives two different arguments. One uses the fact that a semisimple
 Lie algebra is rigid with respect to small deformations; the other uses
 the fact that a semisimple Lie  algebra is its own algebra of derivations,
 so the orbits in the variety of semisimple Lie algebras all have the same
 dimension. 
 
 There is a general procedure for testing when an orbit contains a minimal
 vector, using Hilbert's 1-parameter subgroup criterion for stability in
 the sense of Geometric Invariant Theory \cite{kn:Mum}. In the Lie algebra situation this
 gives a criterion involving the nonexistence of filtrations of a certain
 kind, but the author does not know an easy argument to show that simple
 Lie algebras do not have  such filtrations. However it is also a general fact
 that, in the unstable case, there is a preferred maximally destabilising
 1-parameter subgroup. This theory was developed by Kempf \cite{kn:Ke} and Hesselink \cite{kn:Hes} in the algebraic setting, and in connection with the
moment map and the length function by Kirwan \cite{kn:K} and Ness \cite{kn:N}.
The argument we give in Section 4 is essentially a translation of this theory
into a differential geometric setting. Lauret \cite{kn:L} has studied the application
of this general circle of ideas (Geometric Invariant Theory/Moment maps/minimal vectors) to more sophisticated questions in Lie algebra theory---going beyond
the case of simple algebras. 

One advantage of this method, in the real case, is that the uniqueness of
maximal compact subgroups up to conjugacy emerges as part of the package.
In the usual approach (\cite{kn:Hel}, Theorem 13.5) this is deduced from a
separate argument: Cartan's fixed point theorem for spaces of negative curvature.
We avoid this, although the techniques we apply in Section 4 are very similar in
spirit.

 \section{Appendix}
 
 We give a simple proof of the well-known fact stated in Section 4: that
 the manifold ${\cal H}$ has Property (*). We identify ${\cal H}$ with
 $n\times n$ positive
 definite symmetric matrices of determinant $1$. It suffices to prove the
 statement for the exponential map at the identity matrix. Recall that the
 metric on ${\cal H}$ is given by $\vert \delta H \vert_{H}^{2} = {\rm Tr}\left((
 \delta H)H^{-1}\right)^{2}$. For  fixed
 symmetric matrices $S,\alpha$ and a small real parameter $h$ define
 $$ H(h)= \left(\exp(S+h \alpha) - \exp(S)\right)\exp(-S), $$
 and 
 $$  v= \frac{d H}{dh}\vert_{h=0}.$$
 we need to show that, for any $S$ and $\alpha$, we have
 $$  {\rm Tr}\ v^{2} \geq {\rm Tr}\ \alpha^{2}. $$
 To see this we introduce another real parameter $t$ and set
 $$ H(t,h)= \left( \exp(t(S+h \alpha)) - \exp(tS)\right) \exp(-tS). $$
 Then one readily computes
 $$ \frac{\partial H}{\partial t}= [S, H] + h \alpha \exp( t(S+h \alpha))
 \exp(-tS). $$
 Now differentiate with respect to $h$ and evaluate at $h=0$ to get a matrix
 valued function $V(t)$. Then we have
 $$  \frac{dV}{dt}= \frac{\partial^{2} H}{\partial h\partial t}\vert_{h=0} = [S,V]
 + \alpha. $$
 Clearly $v=V(1)$ and $V(0)=0$, so our result follows from the following
 \begin{lem}
 Let $S,\alpha$ be real, symmetric $n\times n$ matrices and let $V(t)$ be
 the matrix valued function which is the solution of the ODE
 $$    \frac{dV}{dt}= [S, V] +\alpha $$
 with $V(0)=0$. Then  $${\rm Tr}\ V(t)^{2} \geq t^{2}
 {\rm Tr}\ \alpha^{2} $$
 for all $t$.
 \end{lem}
 To see this, consider first a scalar equation
 $$  \frac{df^{+}}{dt\ }= \lambda f^{+} + a, $$
 with $\lambda,a$ constants and with the initial condition $f^{+}(0)=0$. The
 solution is
 $$   f^{+}(t)= \left( \frac{e^{\lambda t}-1}{\lambda} \right) a, $$
 where we understand the expression in brackets is to be interpreted as
 $t$ in the case when $\lambda=0$. Let $f^{-}(t)$ satisfy the similar equation
 $$  \frac{df^{-}}{dt\ } = -\lambda f^{-} + a, $$
 with $f^{-}(0)=0$. Then
 $$   f^{+}(t) f^{-}(t) = t^{2} a^{2} Q(t), $$
 where 
 $$   Q(t)= \frac{(e^{\lambda t} -1)(1- e^{-\lambda t})}{\lambda^{2}t^{2}}=
 \frac{2( \cosh(\lambda t)-1)}{\lambda^{2} t^{2}}. $$
 It is elementary that $Q(t)\geq 1$, so $f^{+}(t) f^{-}(t) \geq t^{2} a^{2}.
 $

  Now consider the operator $\ad_{S}$ acting
 on $n\times n$ matrices. We can suppose $S$ is diagonal with eigenvalues
 $\lambda_{i}$. Then a basis of eigenvectors for $\ad_{S}$ is given by the
 standard elementary matrices $E_{ij}$ and
 $$  \ad_{S}(E_{ij}) = \lambda_{ij} E_{ij}, $$
 where $\lambda_{ij}=\lambda_{i}-\lambda_{j}$.
 Thus the matrix equation reduces to a collection of scalar equations for
 the components $V_{ij}(t)$.
 Since $\lambda_{ji}=-\lambda_{ij}$
 and
 $\alpha_{ij}=\alpha_{ji}$, each  pair $V_{ij}, V_{ji}$ satisfy the conditions
 considered for $f^{+}, f^{-}$ above and we have
 $$  V_{ij}(t) V_{ji}(t) \geq \alpha_{ij}^{2} t^{2}. $$
 (This is also true, with equality, when $i=j$). Now summing over $i,j$ gives
 the result.  
 
 This proof is not very different from the usual discussion of the Jacobi
 equation in a symmetric space. It is also much the same as the proof of
 Helgason's formula for the derivative of the exponential map
 (\cite{kn:Hel}, Theorem 1.7).



\end{document}